\documentclass[10pt,letterpaper]{article}
\usepackage[top=0.85in,left=2.00in,footskip=0.75in,marginparwidth=2in]{geometry}

\usepackage[utf8]{inputenc}

\usepackage{cite}

\usepackage{nameref,hyperref}

\usepackage[right]{lineno}

\usepackage{microtype}
\DisableLigatures[f]{encoding = *, family = * }

\usepackage{changepage}

\usepackage[aboveskip=1pt,labelfont=bf,labelsep=period,singlelinecheck=off]{caption}

\makeatletter
\renewcommand{\@biblabel}[1]{\quad#1.}
\makeatother

\usepackage{lastpage,fancyhdr,graphicx}
\usepackage{epstopdf}
\fancyheadoffset[L]{2.25in}
\fancyfootoffset[L]{2.25in}

\usepackage{color}

\definecolor{Gray}{gray}{.25}

\usepackage{graphicx}

\usepackage{wrapfig}
\usepackage[pscoord]{eso-pic}
\usepackage[fulladjust]{marginnote}
\reversemarginpar

\usepackage{placeins}
\usepackage{parskip}

\newcommand{\asyn}[1]{$\alpha$syn}
\newcommand{\asynuclein}[1]{$\alpha$synuclein}

\begin{document}
\vspace*{0.35in}

\begin{flushleft}
{\Large
\textbf\newline{Mechanistic models of $\alpha$-synuclein homeostasis for Parkinson’s disease: A blueprint for therapeutic intervention}
}
\newline
\\
Elena Righetti\textsuperscript{1,2},
Alice Antonello\textsuperscript{3},
Luca Marchetti\textsuperscript{1,2},
Enrico Domenici\textsuperscript{1,2},
Federico Reali\textsuperscript{1,*}
\\
\bigskip
\bf{1} Fondazione The Microsoft Research - University of Trento Centre for Computational and Systems Biology (COSBI), Rovereto, Italy 
\\
\bf{2} Department of Cellular, Computational and Integrative Biology (CIBIO), University of Trento, Trento, Italy
\\
\bf{3} Department of Mathematics and Geosciences, University of Trieste, Trieste, Italy
\\
\bigskip
* reali@cosbi.eu

\end{flushleft}

\section*{Abstract}
Parkinson's disease (PD) is the second most common neurodegenerative disorder worldwide, yet there is no disease-modifying therapy up to this date. The biological complexity underlying PD 
hampers the investigation of the principal contributors to its pathogenesis.
In this context, mechanistic models grounded in molecular-level knowledge provide virtual labs to uncover the primary events triggering PD onset and progression and suggest promising therapeutic targets. Multiple modeling efforts in PD research have focused on the pathological role of $\alpha$-synuclein (\asyn{}), a presynaptic protein that emerges from the intricate molecular network as a crucial driver of neurodegeneration. Here, we collect the advances in mathematical modeling of \asyn{} homeostasis, focusing on aggregation and degradation pathways, and discussing potential modeling improvements and possible implications in PD therapeutic strategy design.

{\tiny
Keywords: Parkinson’s disease, neurodegeneration, alpha-synuclein aggregation, protein degradation mechanisms, autophagy, mathematical modeling, quantitative systems pharmacology}

\section{Introduction}
Parkinson’s disease (PD) is a debilitating neurodegenerative disorder with currently no disease-modifying therapy or cure \cite{Poewe2017}.
The high prevalence of PD on the global scale and its exponential growth over the last few decades suggest the imminent outbreak of a non-infectious pandemic propelled by aging and industrialization as unconventional vectors of this rapidly growing neurological disorder \cite{Dorsey2018}. Given the enormous strain that a pandemic would put on the healthcare system and the resulting human and economic toll, it is no surprise that investigating PD pathogenesis has rapidly become a serious public health concern. Since 1817, when the “shaking palsy” was first described \cite{Parkinson2002}, considerable progress has been made in understanding PD, from the clinical symptoms to the molecular mechanisms. Such progress has yet to translate into effective treatments that can hamper or even reverse PD progression. 
The current gold-standard therapy (Levodopa or L-Dopa treatment) manages motor symptoms but cannot stop neurodegeneration \cite{Lewitt2016}. It counterbalances dopamine (DA) deficiency in the basal ganglia by relying on the pharmacological input of the DA precursor L-Dopa. Indeed, impairments in the dopaminergic pathway and the subsequent neuronal death in the {\it substantia nigra pars compacta} stand out as primary hallmarks of the disease \cite{Poewe2017}.

\begin{figure}[h!]
	\begin{center}
		\includegraphics[width=\textwidth]{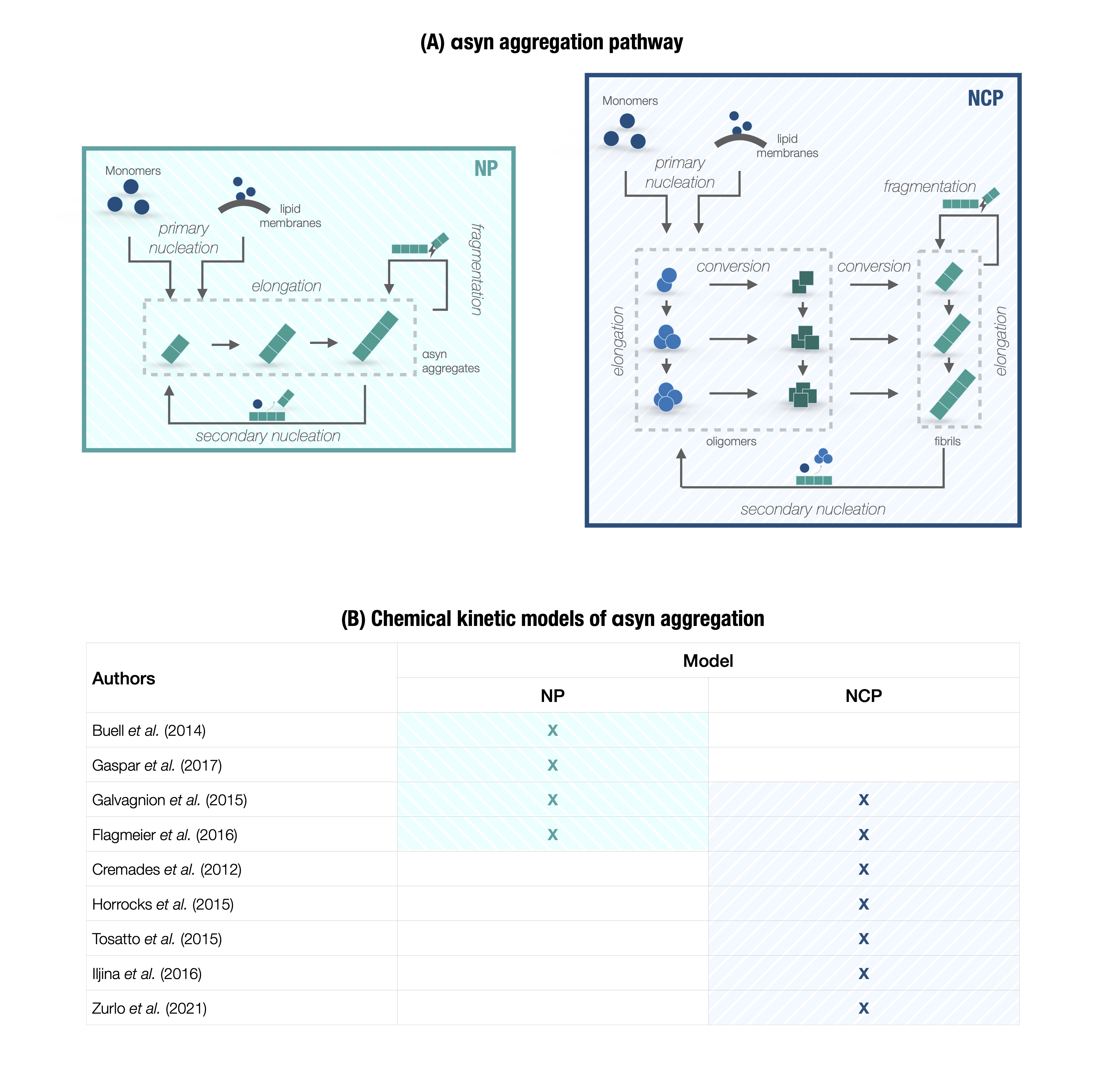}
	\end{center}
	\caption{Chemical kinetic models of \asyn{} aggregation. {\bf (A)} A conceptual scheme of the PD molecular landscape underlying models in Section \ref{Sec2}. The kinetic of \asyn{} aggregation can be represented by either a nucleation-polymerization (NP) process or a nucleation-conversion-polymerization (NCP) process. The pathway includes homogeneous primary nucleation from monomers, heterogeneous primary nucleation catalyzed by lipid membranes, aggregate elongation by monomer addition, structural conversion between oligomeric and fibrillar species (if we consider NCP models), fibril fragmentation, and heterogeneous secondary nucleation on fibril surfaces. Reverse processes, such as monomer dissociation and reverse conversion, are not shown since they are negligible at the early aggregation steps. {\bf (B)} Table collecting the chemical kinetic models of \asyn{} aggregation considered in Section \ref{Sec2} that provide mechanistic insights into PD neurodegeneration. The studies are based on the master equation approach and focus either on NP or NCP models.} 
	\label{CKmodels}
\end{figure}

\begin{figure}
	\begin{center}
		\includegraphics[width=\textwidth]{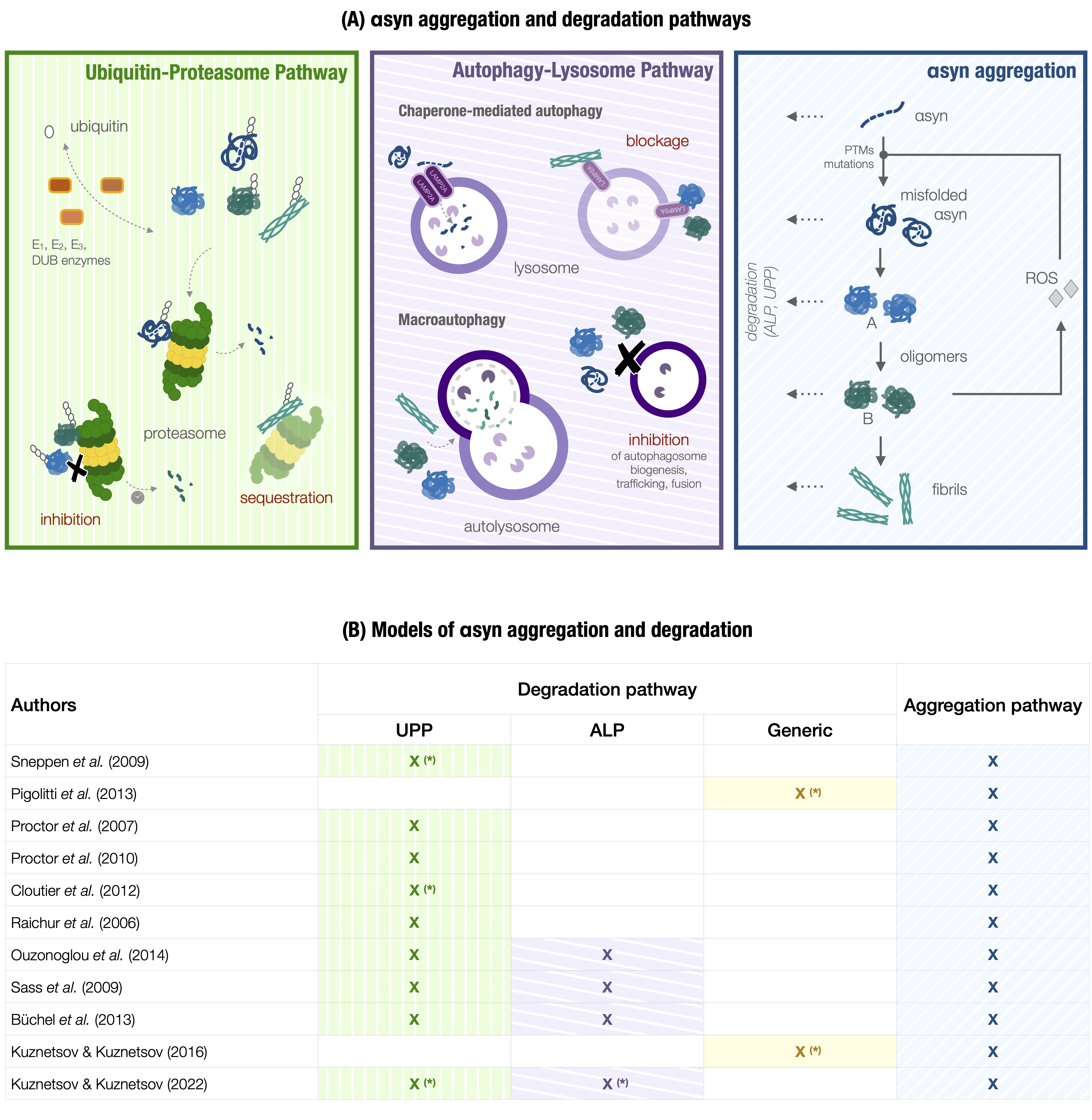}
	\end{center}
	\caption{Models of \asyn{} aggregation and degradation. {\bf (A)} A conceptual scheme of the PD molecular landscape underlying models in Section \ref{Sec3}. A high-level representation of \asyn{} aggregation (light-blue box) is integrated into a more complex network involving the ubiquitin-proteasome pathway (UPP) and the autophagy-lysosome pathway (ALP) (green and purple boxes, respectively), reactive oxygen species (ROS) production, and the corresponding feedback loops. The aggregation process involves \asyn{} misfolding, the conformational change between type-A and type-B oligomers, the conversion from oligomers to fibrils, and degradation of \asyn{} species. Factors triggering \asyn{} misfolding can be post-translational modifications (PTMs), \asyn{} mutations, and high ROS levels. UPP degradation, including ubiquitination and proteasome-mediated degradation, is impaired by fibrils sequestering the proteasome and oligomers inhibiting proteasomal activity either by proteasome overload or directly due to their inherent toxicity. Pathogenic \asyn{} also interferes with ALP components \cite{Xilouri2016}. Fibrils and oligomers can block CMA degradation by obstructing LAMP2A receptors and inhibiting their assembly \cite{Cuervo2004, Martinez2008, Xilouri2009}. Macroautophagy is the main degradation route for \asyn{} aggregates. In return, misfolded \asyn{} and its aggregates can exert an overall inhibitory effect on macroautophagy and hamper its activity by interfering with autophagosome formation, maturation, and trafficking, in addition to lysosomal fusion and transport \cite{Wong2015, Wong2017}. {\bf (B)} Table collecting the models considered in Section \ref{Sec3}.
	The term Generic indicates a generic degradation mechanism, either proteasome or lysosome-mediated. The models denoted by the symbol (*) incorporate degradation mechanisms either as pseudo-first-order reactions or enzymatic reactions; the remaining models are characterized by a more detailed representation of the degradation pathways.} 
	\label{Degmodels}
\end{figure}

Human and molecular genetic research has soon uncovered PD's multifactorial nature, thus unveiling a complex molecular network \cite{Poewe2017, Blauwendraat2020}. In addition to impairments in DA metabolism, numerous processes such as protein aggregation, defective degradation, neuroinflammation, and oxidative stress contribute to disease onset and progression through various feedback mechanisms (see Fujita {\it et al}. \cite{Fujita2014} for a comprehensive disease map). 
The presynaptic amyloidogenic protein $\alpha$-synuclein (\asyn{}) emerges from the intricate molecular landscape as a key driver of PD pathogenesis \cite{Poewe2017}. Here, we focus on two mechanisms governing its homeostasis, {\it i.e.}, the aggregation and degradation pathways. Figures \ref{CKmodels}A and \ref{Degmodels}A display the biological mechanisms involved in these pathways. Their link to PD neurodegeneration is supported by multiple experimental and genome-wide association studies \cite{Blauwendraat2020, Gan2015, Walden2017}. 
Even though multiple studies have highlighted the significant contribution of these processes to PD neurodegeneration \cite{Vekrellis2011, Navarro2020, Olanow2006}, the interplay between the uncontrolled accumulation of \asyn{} aggregates and impairments in the degradation mechanisms is still unclear \cite{Xilouri2013}. 

Given the intricate biological picture outlined so far, identifying the specific role of \asyn{} aggregation and degradation pathways would support therapeutic strategies currently under investigation \cite{Fields2019}, {\it e.g.}, small molecule inhibition, \asyn{} antibodies, and autophagy stimulation, or even pave the way for new pharmacological approaches. 
To this end, the experimental analysis calls for a quantitative systems pharmacology (QSP) approach at the interface between pharmacological research and systems biology \cite{Azer2021}. As recently highlighted in \cite{Geerts2020, Abrams2020, Bloomingdale2021, Bloomingdale2022}, mechanistic QSP models indeed hold a great potential in neuroscience as they provide a systems-level understanding of complex biological processes and insights into their response to drugs.
Specifically, modeling efforts in PD research are moving toward multiscale QSP models that describe molecular, cellular, whole-brain, and organism levels to assess how the effects of a molecular perturbation can scale up to influence the clinical outcome \cite{Stephenson2022}. 

Starting from the molecular level, various mechanistic models have tackled PD complexity. 
As highlighted by previous studies on mechanistic models for neurodegenerative disorders \cite{LloretVillas2017, Bakshi2019}, they have mainly focused on protein aggregation, synaptic transmission, apoptosis, oxidative stress, and genetic components rather than protein degradation. The relatively small number of models including degradation processes is likely to rise along with increasing evidence that supports the critical role of the protein clearance machinery in PD \cite{Navarro2020, Klein2018}.
With this in mind, we collect the advances in mathematical modeling of \asyn{} aggregation and degradation by the ubiquitin-proteasome pathway, macroautophagy, and chaperone-mediated autophagy. 
We divide the mechanistic models into two groups: (i) single-pathway models of \asyn{} aggregation representing the chemical kinetics of the process ({\it i.e.}, {\bf chemical kinetic models});
(ii) multiple-pathway models of \asyn{} homeostasis focusing on the interplay between \asyn{} aggregation and degradation mechanisms ({\it i.e.}, {\bf degradation models}).
We discuss their potential improvements and possible implications in PD therapy development. 


\FloatBarrier

\FloatBarrier


\section{Mechanistic models of \texorpdfstring{\asyn{}}{} aggregation kinetics}
\label{Sec2}

Protein aggregation and misfolding are central features in the molecular network of neurodegeneration \cite{Chiti2017}. It is therefore no surprise that a large part of the modeling effort targeting PD investigates the chain of microscopic events in \asyn{} aggregation. These studies fall into a broader line of research devoted to uncovering the generality of amyloid formation by focusing on the chemical and physical aspects of this phenomenon \cite{Dobson2003, Knowles2009}.
Here, we analyze mechanistic models that extend the formalism of chemical kinetics to the aggregation process; we thus refer to them as {\bf chemical kinetic models} (see Figure \ref{CKmodels}). These models adopt a master equation approach building on the pioneering work by Oosawa $\&$ Kasai \cite{OosawaKasai1962} to provide a mathematical formulation of each molecular step in the process and derive the corresponding rate constants.
To gain novel biological insights and formulate new scientific queries, {\it in silico} models evolve along with increasingly advanced experimental procedures in a virtuous cycle of model calibration and validation against {\it in vitro} quantitative measurements \cite{Chiti2017}.

\subsection{\texorpdfstring{\asyn{}}{}  misfolding and aggregation and PD neurodegeneration}

Soluble \asyn{} is an aggregation-prone protein with an intrinsically disordered conformation that prevents its spontaneous self-assembly \cite{Gracia2020}. Many factors can trigger \asyn{} misfolding and accumulation, such as high \asyn{} concentrations, genetic mutations, increased reactive oxygen species (ROS) levels, and interaction with lipid membranes \cite{Hijaz2020}. 
The process initiates upon structural changes that expose the highly amyloidogenic region of the protein to the cytosol. Misfolded monomers slowly aggregate into oligomeric species, thus creating nucleation seeds. In return, seeding-competent oligomers can rapidly accumulate and elongate by monomer addition forming protofibrils and insoluble fibrils (Figure \ref{CKmodels}A). Multiple \asyn{} species eventually gather into amyloid inclusions known as Lewy bodies (LBs) and Lewy neurites (LNs), which represent one of PD's main molecular hallmarks \cite{Spillantini1997}. 

A considerable body of research points to the harmful effect of \asyn{} aggregation and its crucial contribution to PD neurodegeneration. For instance, the {\it SNCA} gene, encoding for \asyn{}, was the first gene associated with PD \cite{Polymeropoulos1997}. 
Still, there is no agreement on the specific mechanisms of toxicity. 
Seeding and spreading ability and inherent toxicity are the functional and structural factors that determine the overall toxicity of aggregates \cite{Gracia2020, Froula2019}. 
According to {\it in vitro} and {\it in vivo} experimental evidence, the most cytotoxic species are oligomers with a partially folded structure similar to mature fibrils \cite{Hijaz2020, Froula2019, Cremades2017}. Due to their conformation, some oligomeric species can create pores within the cell membrane, thus enhancing its permeability and dysregulating calcium homeostasis \cite{Alam2019}. They can also impair mitochondrial function and proteolytic machinery. On the other hand, small fibril fragments have been identified as major-spreading agents in PD because of their high seeding efficiency and internalization ability \cite{Hijaz2020}. 
According to one of the leading hypotheses for \asyn{} propagation, pathogenic species are assumed to spread in a prion-like manner \cite{Brundin2010}; they are released in the extracellular space and then endocytosed by neighboring neurons, where they recruit non-native monomers acting as seeds for further aggregation. As a result, protein aggregates propagate across interconnected brain areas over decades, eventually resulting in PD neurodegeneration \cite{Poewe2017}.

Considerable progress has been made in understanding the kinetic of \asyn{} misfolding and aggregation, however, there is no agreement on many aspects of this mechanism. Chemical kinetic models can provide valuable insights in understanding, for example, the relative contribution of different microscopic events to the process, the structural and functional features of \asyn{} aggregate species \cite{Cremades2017}, and the relationship between \asyn{} aggregation and LB formation and maturation \cite{Hardenberg2021, Shahmoradian2019}.

\subsection{The master equation approach}
The aggregation pathway of a generic amyloidogenic protein can be represented by a minimal set of fundamental reactions, regulating the transition from soluble monomers to fibrils through various intermediate states. As shown in Figure \ref{CKmodels}A, we can differentiate between primary reactions, {\it i.e.}, homogeneous primary nucleation and elongation, and secondary reactions, {\it i.e.}, heterogeneous secondary nucleation and fibril fragmentation \cite{Cohen2011a, Meisl2016}. Homogeneous nucleation initiates the aggregation process by combining monomers into newly formed aggregates that, in turn, elongate by monomer addition. By integrating these reactions, fibril fragmentation and heterogeneous nucleation on aggregate surfaces multiply the aggregate number in a positive feedback loop with elongation, thus ensuring self-replication. Their relative contribution to the aggregation process depends on the specific protein considered \cite{Knowles2009}. Primary and secondary reactions can incorporate more complex events, such as interaction with lipid membranes and multistep nucleation.
The master equation approach \cite{OosawaKasai1962, Michaels2014} translates the reaction network into an infinite system of coupled non-linear differential equations describing, in probabilistic terms, the time course of the species involved, from monomers to filaments identified by their monomer count.

The main goal of this approach is to find an analytical solution for the system to describe its reaction time course. As a result, the analytical solution can be fitted to experimental data to identify kinetic rate constants and provide mechanistic insights into the aggregation process. Since the master equation formulation is often analytically and numerically intractable, specific approximation methods can be employed to reduce the system to two differential equations describing the time evolution of aggregate number and mass concentrations, {\it i.e.}, readily available observables. 
Multiple efforts have channeled into a formal analysis of this minimal closed system. 
Knowles and colleagues \cite{Knowles2009, Cohen2011a, Cohen2011b} proposed fixed-point approaches to derive explicit solutions for the entire time course of aggregates. All the identified parameters had a clear physical meaning related to the underlying molecular events. Furthermore, analytical solutions exhibited scaling laws used in model fitting to determine the importance of primary and secondary reactions \cite{Michaels2016}. For instance, in the case of breakable fibrils \cite{Knowles2009}, the authors found that the reaction time course mainly depended on a parameter combining elongation and fragmentation rate constants. By analyzing the relationship between this parameter and the observables, such as half-time and maximal growth rate, the authors showed that even small perturbations in the fragmentation reaction significantly affected the aggregation process. This procedure was implemented in a global fitting tool for chemical kinetic models of {\it in vitro} protein aggregation (AmyloFit) \cite{Meisl2016}. 

Model parsimony and data accuracy are essential requirements for the fitting protocol. A parsimonious model can be fitted to a wide range of experimental data acquired for different initial monomer concentrations and in the presence and absence of preformed aggregate seeds. Techniques for measuring aggregate species vary from fluorescent dyes ({\it e.g.}, thioflavin T) to mass spectrometry for aggregate mass concentration and single-molecule fluorescence for oligomer concentrations \cite{Meisl2016, Cremades2012, Perrett2014}. With improvements in experimental approaches come recent advances in the mathematical framework, including the investigation of stochastic aggregation kinetics \cite{Michaels2018}, the spatial dependence of protein aggregation \cite{Cohen2014}, and the nature of oligomeric species \cite{Cremades2017, Dear2020a}.

\subsection{Chemical kinetic models of \texorpdfstring{\asyn{}}{}  aggregation}

A significant number of chemical kinetic models address \asyn{} aggregation and its implications in PD. 

{\it Nucleation-polymerization (NP) models.} NP models describe aggregation as fibril formation by monomer addition, which follows a nucleation event and escalates by secondary reactions. This process is represented by a characteristic sigmoidal growth curve, identified by a lag phase, a growth-phase, and a plateau \cite{Knowles2009}. 
By relying on the kinetic theory of aggregation that revolves around the master equation approach, multiple studies have combined experimental analysis with these mechanistic models to investigate underlying microscopic events and the effect of varying experimental conditions. 

Buell {\it et al}. \cite{Buell2014} focused on the effect of pH on \asyn{} aggregation in the presence and absence of preformed fibril seeds. Homogeneous and heterogeneous nucleations and fragmentation were negligible in purely monomeric solutions at neutral pH. Adding preexisting aggregates to the solution enhanced the aggregation rate and shortened the lag phase of the sigmoidal growth curve. To obtain the elongation rate constant, the authors used a minimal model involving only an elongation reaction saturating at high monomer concentrations ({\it i.e.}, two-step elongation). On the other hand, a mildly acidic pH uncovered a much faster aggregation, that is, a longer lag phase and a more rapid growth phase, due to either surface-catalyzed nucleation or fibril fragmentation.
Building on \cite{Buell2014}, Gaspar {\it et al}. \cite{Gaspar2017} determined the pH-dependent role of these two events by relying on the theoretical analysis of a model including primary and secondary reactions. According to an experimental analysis supported by a global fitting procedure \cite{Meisl2016}, the autocatalytic nature of aggregation was mainly determined by secondary nucleation of monomers on fibril surfaces. These studies provided insightful information on physiological mechanisms; for example, they showed that even a slight pH change associated with different intracellular environments may enhance the production of toxic aggregates. 

Small lipid vesicles can alter \asyn{} aggregation as much as pH levels \cite{Galvagnion2015}. Upon interaction with lipid membranes, the protein adopts a stable $\alpha$-helical conformation \cite{Fusco2014}. Galvagnion {\it et al.} \cite{Galvagnion2015} showed that an imbalance between free and membrane-bound \asyn{} states can trigger amyloid formation by enhancing heterogeneous primary nucleation on lipid surfaces. To confirm and quantify these results, the authors employed an NP model including this nucleation reaction, two-step elongation, and the conversion between intermediate species \cite{Cremades2012}. 
Relying on \cite{Buell2014, Galvagnion2015}, Flagmeier {\it et al}. \cite{Flagmeier2016} investigated how mutations in \asyn{} influence the fate of amyloid fibrillation. The authors performed three experiments to separately analyze initiation, elongation, and multiplication steps for all \asyn{} variants. Each experiment was associated with a mechanistic model: the NP model from \cite{Galvagnion2015} for lipid-induced aggregation; and linear polymerization models for strongly and weakly seeded aggregation processes at neutral and mildly acidic pH. Experimental and theoretical results agreed that the rate of surface-catalyzed nucleation on lipid membranes and fibril surfaces could vary by several orders of magnitude across different mutations. 

{\it Nucleation-conversion-polymerization (NCP) models.} The NP models are coarse-grained to the extent that they do not explicitly account for \asyn{} intermediate species. The modeling choice mainly relates to the absence of accurate measurements of oligomer concentrations in the aggregation process. Oligomers are indeed transient and highly heterogeneous species. The development of new methods to report their levels has motivated the transition to NCP models \cite{Dear2020a, Garcia2014}. These models include the structural interconversion from disordered clusters to oligomers with a partially formed fibrillar structure, eventually growing into fibrils.

Multiple studies on \asyn{} aggregation have provided experimental and theoretical support to this updated kinetic theory \cite{Horrocks2015, Tosatto2015, Iljina2016, Zurlo2021}. 
Single-molecule F{\"o}rster Resonance Energy Transfer experiments can distinguish two oligomeric populations associated with low and high FRET levels. Type-A oligomers are early-formed, disordered, highly degradable, and mildly toxic species, whereas type-B oligomers are later-forming, more compact, degradation-resistant, and ROS promoting species \cite{Cremades2012, Iljina2016}. Both oligomeric species appear as required subsequent steps in the aggregation process \cite{Dear2020b}. 
As displayed in Figure \ref{CKmodels}A, NCP models tailored for \asyn{} aggregation include primary nucleation, which generates type-A oligomers from monomers, and the structural conversion from type-A to type-B oligomers and then fibrils \cite{Cremades2012, Iljina2016}. Both fibrillar and oligomeric species increase in size by monomer-dependent elongation.
In this reaction network, the conformational change between oligomeric species emerges as a crucial step for initiating \asyn{} aggregation. Upon fitting the model to experimental data, Cremades {\it et al}. \cite{Cremades2012} showed that the corresponding kinetic rate is slow, thus providing a sufficient lag time for the protein quality control machinery to hamper the pathological process. A potential therapeutic strategy may employ molecular chaperones or small molecules to regulate this reaction. A follow-up study \cite{Iljina2016} used the NCP model to quantify seeding ability and conditions for type-B oligomers and fibrils, with implications on prion-like spreading. Model predictions showed that fibrils were more prone to seed aggregation, whereas oligomers contributed mostly to the increase of oxidative stress.

The studies reported above focused on the early stages of \asyn{} aggregation; thus, they did not require full-time analytical solutions. To extend the analysis to the entire time course, Dear {\it et al.} \cite{Dear2020a} provided a general mathematical framework for NCP models. The authors proposed a reaction network with a single oligomeric state.
Such a coarse-grained model is suitable to avoid overfitting in the analysis of experimental data on total oligomer concentration that, as such, did not account for differences between oligomeric populations, {\it e.g.}, see \cite{Zurlo2021}. 
Moreover, the model is instrumental for classifying oligomers according to theoretically defined metrics: persistence or half-life, abundance, and productivity, {\it i.e.}, the propensity to form fibrils rather than dissociate. The authors proved that oligomers were more prone to dissociation than fibrillation as a universal feature. Focusing on \asyn{} oligomers, results showed that they were very kinetically stable, more persistent than monomers, relatively abundant, and highly productive when compared to other proteins. 

Overall, chemical kinetic models contribute to a detailed mechanistic understanding of \asyn{} aggregation by combining theoretical results and {\it in vitro} data acquired through advanced techniques. By exploring a wide range of {\it in silico} experimental conditions, these models return coherent predictions that can be useful for studying PD. Further investigation in this direction should focus on whether and how these results apply to {\it in vivo} settings to provide insights into the mechanisms of neurodegeneration and, possibly, inform therapeutic strategy design \cite{Sinnige2022}.

\section{Mechanistic models of \texorpdfstring{\asyn{}}{} aggregation and degradation pathways}
\label{Sec3}

Given the central role of degradation in regulating \asyn{} homeostasis, multiple models have investigated this process in relation to PD neurodegeneration, including models solely focused on \asyn{} aggregation and degradation or {\it miscellaneous} models that include additional pathogenic processes ({\it e.g.}, ROS production and defective dopamine metabolism). 
These {\bf degradation models} offer a heterogeneous picture in terms of specific mechanisms they consider, granularity, and modeling approaches \cite{Bakshi2019} (see Figure \ref{Degmodels}). The biological processes can be represented by varying degrees of detail, and the modeling formalism can range from ordinary differential equations (ODEs) and algebraic equations to stochastic simulations, partial differential equations, or flux balance analysis \cite{LloretVillas2017,Bakshi2019}. 

\subsection{Defective protein degradation machinery and PD neurodegeneration}

The proteolytic machinery combined with the molecular chaperone system ensures the renewal of functional components and regulates the fate of altered proteins to avoid the build-up of intracellular damage. 
In PD and other neurodegenerative diseases, degradation efficiency declines with age and oxidative stress \cite{Szweda2002}, thus increasing misfolded \asyn{} and aggregate levels.
On the other hand, \asyn{} aggregation can hamper the proteolytic machinery causing the accumulation of other cellular components \cite{Tai2008}. The interplay of \asyn{} aggregation and impairment of degradation in multiple negative feedbacks poses a causality dilemma \cite{Ciechanover2003}: which pathological process comes first? There is still no clear answer. Degradation defects and \asyn{} aggregation can be either primary or secondary events of neurodegeneration \cite{Xilouri2013,Scrivo2018}. The investigation of their role in PD is even more challenging due to interfering factors such as increased oxidative stress and mitochondrial dysfunctions \cite{Poewe2017, Szweda2002}. QSP models can help in elucidating these contributions by providing {\it in silico} benchmarks to assess these hypotheses. 

{\it The ubiquitin-proteasome pathway (UPP).} UPP is the dominant clearance route for mislocated, misfolded, damaged, and short-lived proteins. Target substrates are selectively tagged by polyubiquitin chains through multiple enzymatic reactions and then transferred to the proteasome. Ubiquitinating and deubiquitinating enzymes such as PARKIN and UCH-L1 play a central role in modulating protein degradation in PD, as indicated by UPP-related familial cases \cite{Walden2017, Tai2008}. 
For instance, UCH-L1 maintains a stable pool of monomeric ubiquitin; an imbalance of this pool can lead to poor substrate labeling with consequent impaired proteasomal degradation.
Furthermore, the presence of LBs as a PD molecular hallmark corroborates the correlation between proteasomal dysfunctions and PD pathogenesis. Indeed, LBs mainly consist of ubiquitin-tagged \asyn{} aggregate, ubiquitinated proteins, and UPP components 
\cite{Tai2008, Kuzuhara1988, Shults2006}. According to multiple studies \cite{Olanow2006, McNaught2001}, their composition suggests a potential cytoprotective role for LBs, which may act as sinks for poorly degraded proteins to reduce the toxic effect of protein accumulation and slow down neurodegeneration. Still, it is unclear whether these cytoplasmic inclusions represent a harmful or defense molecular mechanism for the neuron.
UPP dysfunctions and \asyn{} aggregation create a negative feedback loop. Genetic and age-related impairments in UPP may hamper the proper degradation of misfolded aggregate-prone \asyn{}. On the other hand, \asyn{} species refractory to proteasomal degradation, such as mutant \asyn{} forms, oligomers, and fibrils, can inhibit proteasomal activity by sequestering the proteasome. In addition, the interaction with oligomers and oxidatively damaged \asyn{} may directly alter proteasome subunits, as shown in the green box in Figure \ref{Degmodels}A.

{\it The autophagy-lysosome pathway (ALP).} ALP removes long-lived proteins, lipids, and dysfunctional organelles, such as aging mitochondria ({\it i.e.}, mitophagy), thus maintaining neuronal homeostasis. It includes three degradation mechanisms: macroautophagy, chaperone-mediated autophagy (CMA), and microautophagy. 
Experimental findings indicate that both UPP and ALP degradation of \asyn{} can occur \cite{Webb2003}. The autophagic mechanisms are involved in physiological \asyn{} turnover and \asyn{} aggregates degradation (see the purple box in Figure \ref{Degmodels}A). Specifically, wild-type \asyn{} monomers are degraded mainly by CMA \cite{Cuervo2004}, which relies on molecular chaperones to recognize a specific motif on the target substrate. Upon chaperone binding, the cargo interacts with the membrane receptor LAMP2A, which assembles into multimeric translocation complexes. The substrate then enters the lysosome lumen, where it is rapidly degraded. CMA also clears ALP and UPP-related proteins, such as LRRK2, UCH-L1, DJ-1, and GCase \cite{Hou2020}. 
Furthermore, macroautophagy represents the main clearance route for \asyn{} aggregates due to the inaccessibility of the CMA degradation motif \cite{Scrivo2018}. This degradation mechanism consists in the formation of autophagosomes, {\it i.e.}, double-membrane vesicles surrounding the cargo and fusing with lysosomes to allow hydrolase-mediated degradation in the lysosomal lumen. In selective macroautophagy, specific autophagosome membrane receptors target protein aggregates. 
The link between impaired autophagy and PD neurodegeneration conceals a complex interaction network between \asyn{} aggregates and the lysosomal proteolytic machinery. Growing evidence supports the negative effect of \asyn{} aggregation on autophagy acting on multiple levels \cite{Scrivo2018, Xilouri2016}, as displayed in Figure \ref{Degmodels}A. 
At the same time, defective ALP cannot handle harmful \asyn{} species, thus leading to their accumulation. Such defects can originate from various sources, such as increased oxidative stress, aging, and dysfunctional proteins related to genetic PD forms, namely, LRRK2, DJ-1, and GBA \cite{Blauwendraat2020}, and with a physiological role in ALP.

\subsection{Degradation models}

{\it UPP models.} 
Relying on experimental evidence pointing to the essential contribution of UPP to neurodegeneration \cite{McNaught2004, Cuervo2010}, Sneppen {\it et al}. \cite{Sneppen2009} provided a minimal theoretical model of the interplay between impairments in proteasomal degradation and \asyn{} fibrillation. The proteasome ensured protein homeostasis by keeping the amount of \asyn{} aggregates low. On the other hand, insoluble \asyn{} fibrils sequestered this proteolytic machinery and formed a slowly-degraded complex, thus reducing proteasome function. 
Computational simulations of this ODE system supported by stability analysis showed a bifurcation behavior. When the fibril level exceeded a particular threshold related to proteasome production capacity, the double-negative feedback broke \asyn{} homeostasis in favor of a pulsatile oscillatory regime characterized by spikes and long periods of low proteasome concentrations. As a result, the system alternated between the hypothetical recovery of proteolytic activity and the uncontrolled accumulation of harmful aggregates and fibrils. Several oscillations eventually built up toxicity over decades, thus explaining the slow progression of PD neurodegeneration. The authors therefore associated proteasome overload and the subsequent oscillatory regime with the onset of sporadic PD phenotype. 
In a follow-up study \cite{Pigolotti2013}, the stochastic analysis of the core negative feedback confirmed and generalized the above phenomenon for the whole protein clearance machinery. The bifurcation behavior characterized by a {\it healthy} homeostatic state and a {\it disease} oscillatory state indeed emerged for both UPP and ALP degradation. In addition, stochastic simulations captured dynamical differences between the two degradation systems, modeled as enzymatic reactions slowly reducing aggregate concentrations. Lysosomes were assumed to be more efficient but less abundant than proteasomes. Fibril-mediated sequestration of the proteasome therefore led to regular and almost deterministic oscillations, whereas lysosome overload resulted in irregular fluctuations.

In contrast to the high-level representation of the UPP in \cite{Sneppen2009, Pigolotti2013}, Proctor {\it et al.} \cite{Proctor2007, Proctor2010} proposed a detailed description of this pathway by including mono- and polyubiquitination, proteasome binding, ATP-dependent degradation, and proteasome inhibition by protein aggregates. These studies were devoted to the mathematical modeling of the aging process, which contributes to the decline of proteolytic activity and the build-up of oxidative stress due to flawed mitochondria  \cite{Kirkwood2003, Mc2017}. Both models indeed investigated the role of the proteasomal degradation pathway in age-related neurodegeneration. 
In \cite{Proctor2007}, Proctor and colleagues first presented a stochastic model for aggregation and proteasomal degradation of a generic protein under normal homeostasis and aging. The authors then tailored the system to assess the UPP role in \asyn{} aggregation during PD neurodegeneration \cite{Proctor2010}. The final model included additional modules representing \asyn{} metabolism and UCH-L1 deubiquitination activity in proteasomal degradation. This stepwise approach relied on the idea that the molecular landscape of neurodegeneration can be arranged in building blocks, as most biological processes are assumed to be modular \cite{LloretVillas2017, Hartwell1999}.
Both general and specific models proved to be functional for investigating the negative feedback of protein aggregate formation and UPP impairment. Employing the more general model \cite{Proctor2007}, the authors explored the impact of various experimental procedures, {\it e.g.}, shutting down a ubiquitin-protein ligase or the proteasome activity. The results suggested that depleted ubiquitin pools might contribute to proteostasis impairment. In addition, relying on the model tailored for PD \cite{Proctor2010}, the authors analyzed the effects of both mutated and oxidatively damaged UCH-L1. Stochastic simulations indicated that proteasome inhibition might not be the triggering event for neurodegeneration but rather an exacerbating factor accelerating \asyn{} aggregation and inclusion formation. 
Both models were validated against data gathered in {\it ad-hoc} designed experiments, in a virtuous cycle of model refinement and {\it in vitro} experiments.

{\it UPP and ROS models.}
Since oxidative stress and impairments in UPP degradation are interconnected, multiple studies have analyzed their synergistic effect on \asyn{} homeostasis.
Cloutier {\it et al.} \cite{Cloutier2012a} provided a mathematical model accounting for the interplay between \asyn{} aggregation, ROS production, and proteasome-mediated degradation. Given the variety of molecular processes involved, the authors were able to evaluate the impact of multiple risk factors on \asyn{} metabolism. For instance, perturbations of system parameters related to oxidative stress and protein clearance capacity reproduced the high variability of LB concentration in PD forms \cite{Cloutier2012a, Foltynie2002}.
Moreover, supported by experimental verification \cite{Finnerty2013}, deterministic simulations of the model uncovered a bistable process whereby the system irreversibly switched from a healthy to a disease state, identified by high misfolded \asyn{} and ROS levels, in response to toxin exposure, genetic mutations, and aging. On the other hand, the model could not recapitulate the oscillatory regime due to decreased levels of available clearance machinery obtained in \cite{Sneppen2009, Pigolotti2013}. 
To identify the primary mechanisms determining the bifurcation behavior, the authors reduced the system to focus only on the double-positive feedback between \asyn{} misfolding and ROS production \cite{Cloutier2012b}. This minimal version retained the bistable behavior, thus associating PD onset with a switch-like transition in response to enhanced ROS concentrations. In line with \cite{Sneppen2009, Pigolotti2013}, the central idea remained that neurodegeneration results from the transition between two dynamical behaviors. 

Raichur {\it et al}. \cite{Raichur2006} also presented a model targeting the interplay between oxidative stress, \asyn{} aggregation, and UPP dysfunctions.
Compared to Cloutier {\it et al}. \cite{Cloutier2012a}, this system displayed a more granular description of the UPP pathway, from ubiquitination and proteasome recognition to degradation. 
The close molecular detail of ROS production and UPP degradation provided in the model inevitably hampered the theoretical analysis. 
On the other hand, the explicit representation of proteasome-mediated degradation enabled the authors to explore different scenarios of UPP-related genetic predisposition: mutations in {\it UCH-L1}, {\it PARKIN}, and {\it DJ-1} genes, combined with age-related oxidative stress.
These familial PD forms determined the accumulation of ubiquitinated proteins, leading to a more rapid \asyn{} aggregation than sporadic forms. Furthermore, the synergistic combination of ROS production and UPP-related genetic mutations fueled \asyn{} aggregation. These results pointed to therapeutic strategies involving restored UPP functionality or enhanced ALP activity in combination with antioxidant treatments.

{\it UPP and ALP models.}
The modeling picture outlined so far focuses on the ubiquitin-proteasome pathway, while the autophagy-lysosome pathway is under-represented despite its critical role in \asyn{} turnover; when considered, 
autophagic clearance has often been incorporated as a pseudo-first-order reaction.
This is not the case of Ouzounoglou {\it et al}. \cite{Ouzounoglou2014} who provided a detailed representation of selective macroautophagy, CMA, and proteasome-mediated degradation, focusing on the inhibitory feedback between \asyn{} oligomers and the CMA receptor LAMP2A.
Upon model calibration and validation against \asyn{} overexpression data from human neuroblastoma cells, the authors performed stochastic simulations to explore different interventions with neuroprotective potentials, such as shutting off DA production, increasing LAMP2A levels, and reducing \asyn{} synthesis. The resulting predictions of increased cell viability were in agreement with experimental evidence and literature \cite{Xilouri2009, Vekrellis2009}.

To provide valuable insights into the role of \asyn{} degradation in neurodegeneration, a mechanistic model should integrate ALP and UPP degradation mechanisms, ideally accounting for the crosstalk between these pathways.
In this direction, Sass {\it et al}. \cite{Sass2009} proposed a model targeting asyn{} aggregation, UPP and ALP degradation, and DA metabolism and based on the biochemical system theory. A granular representation of these mechanisms and their interactions facilitated the analysis of specific system perturbations in both disease and treated states ({\it e.g.}, enhanced \asyn{} fibrillation and aggregation rate). Including DA metabolism in the model also enabled the authors to simulate impairment of DA vesicle packaging and available treatments related to its shortage, such as L-Dopa therapy combined with monoamine oxidase inhibitors. 
Increasing model complexity even further, B{\"u}chel {\it et al}. \cite{Buchel2013} presented an extended model of a dopaminergic neuron including the highest number of molecular components among the available PD models of \asyn{} homeostasis \cite{LloretVillas2017}. The model described eleven biological processes, including DA metabolism and transport, \asyn{} aggregation, oxidative stress, UPP and ALP clearance mechanisms, and specific details on PARKIN and DJ-1 proteins. Here, the modeling approach adopted was crucial to tackle model complexity. The authors employed flux balance analysis to investigate the steady-state behavior of the model rather than deterministic or stochastic analysis based on mass action kinetics. The qualitative prediction of multiple {\it in silico} experiments pointed to increased ROS, neurotoxin MPTP, and \asyn{} levels as crucial factors in PD neurodegeneration.

{\it \asyn{} degradation and axonal transport models.}
The degradation models analyzed so far do not account for the spatial dependence of \asyn{} aggregation and degradation mechanisms. However, \asyn{} homeostasis and disease progression may strongly depend on intraneuronal protein propagation. Along this line, Kuznetsov $\&$ Kuznetsov \cite{Kuznetsov2016a,Kuznetsov2016b} investigated \asyn{} axonal transport and how its perturbation between the somatic and synaptic compartments affected the aggregation and degradation processes. By relying on a two-compartment ODE model \cite{Kuznetsov2016b}, the authors predicted axonal death due to \asyn{} accumulation at the synaptic terminal, pointing to defective protein clearance machinery as the primary cause for \asyn{} aggregation onset compared to reduced \asyn{} production and transport \cite{Kuznetsov2016b}. Differences in the degradation of \asyn{} monomers and aggregates are determined by assuming distinct half-lives.
Furthermore, building on recent findings on the absence of \asyn{} fibrils in LB inclusions \cite{Shahmoradian2019}, the authors proposed a mathematical model of \asyn{} aggregation, degradation, and transport to test the hypothetical autocatalytic nature of LB formation and its impact on protein accumulation \cite{Kuznetsov2022}. Here, the authors assumed spatial-dependent proteasomal and autophagic degradation mechanisms. The study pointed to the central role of the autophagic clearance pathway in LB formation and PD progression; defective autophagy may favor the aggregation of membrane-bound organelles catalyzed by \asyn{} aggregates and, thus, LB formation. Most degradation models do not differentiate between \asyn{} fibrillation and LB formation, whereas chemical kinetic models do not even include these cytoplasmic inclusions. The detailed implementation of LB formation and maturation would allow for the investigation of LB neuroprotective role \cite{McNaught2001, Kuznetsov2022}.

The wide variety of degradation models reflects the multifactorial nature of PD. Such a modeling effort has provided multiple insights into the interplay between \asyn{} aggregation and protein degradation by spanning various levels of molecular detail, from coarse-grained to highly granular representations. 
Minimal models have been used to uncover dynamical features of disease progression ({\it e.g.}, \cite{Sneppen2009, Pigolotti2013, Cloutier2012b}). More detailed models ({\it e.g.}, \cite{Proctor2010, Ouzounoglou2014}) have proven to be suitable for testing the effects of potential treatments. Overall, degradation models move toward an integrative representation of PD pathogenesis. Future work may focus on the crosstalk between UPP and ALP pathways, the interplay between autophagy and apoptosis in PD neurodegeneration \cite{Tavassoly2015}, and the integration of other control systems for protein homeostasis such as the molecular chaperone system \cite{Proctor2005}.  


\section{Discussion}
Mechanistic modeling grounded in molecular-level knowledge is a powerful tool for tackling PD biological complexity. By exploring the intricate interaction network of neurodegeneration through {\it in silico} experiments, mathematical models can investigate the key drivers of PD pathogenesis, discover promising targets, and even simulate the effect of potential single or combination therapies.
Given the pivotal role of \asyn{} homeostasis in PD, this review has provided a collection of available mathematical models of \asyn{} aggregation and degradation pathways. 

Chemical kinetic models focus on a single molecular pathway and rely on a unifying kinetic theory that translates into the consistency of modeling formalism ({\it i.e.}, master equations approximated to ODE systems). In contrast, degradation models describe multiple biological processes and adopt a wide range of approaches to analyze the dynamical features of the system. 
The latter are modeled as deterministic and stochastic simulations based on mass action kinetics; yet, there are some exceptions, such as flux balance analysis, mainly related to the complexity of the reaction network considered. Stochastic approaches \cite{Gillespie1977, Wilkinson2018} have been employed to analyze biological processes with low molecular counts and capture their inherent fluctuations.
Future investigation may explore the role of intrinsic noise in PD onset and progression since it can unpredictably affect the system behavior, for example, by triggering the transition from a healthy to a disease state leading to PD onset \cite{Bakshi2019, Cloutier2012a}.

{\it Integrating chemical kinetic and degradation models}.
Both chemical kinetic and degradation models include a representation of \asyn{} aggregation. As shown in Figures \ref{CKmodels}A and \ref{Degmodels}A (light-blue box), chemical kinetic models account for each molecular step separately, whereas degradation models employ macro-variables for different \asyn{} aggregate species according to their size and functional features. The master equation approach returns parameters that physically relate microscopic events to macroscopic observables. On the other hand, a higher-level representation of the process enables the targeting of additional pathways and the investigation of their interplay in the molecular network. A modeling opportunity lies at the intersection of the two groups: degradation models can be integrated with a chemical kinetic representation of the aggregation process. 
In this way, they can benefit from the kinetic rate constants derived from fitting to the wide variety of {\it in vitro} data available. At the same time, chemical kinetic models can be extended to include UPP and ALP pathways, converging toward an {\it in vivo} representation.
These integrated models can be employed to determine the relative importance of each microscopic event in relation to degradation mechanisms ({\it e.g}., the structural interconversion between type-A and type-B oligomers), thus suggesting specific target species and reactions. 
Recent studies have moved a few steps toward model integration \cite{Pigolotti2013, Thompson2021}.
Here, clearance mechanisms are approximated by pseudo-first-order reactions to facilitate theoretical analysis. Further investigation should focus on an explicit representation of the degradation processes, especially the ALP pathway, considering the growing interest in the role of autophagy in neurodegeneration \cite{Gan2015, Klein2018, Senkevich2020}. However, this choice translates into an increased complexity of the system itself. 
Future work should favor an integrative representation of pathological mechanisms that balances high granularity and analytical or computational tractability \cite{Ribba2017}.

{\it In vitro/in vivo translation}.
Integrating chemical kinetic and degradation models is part of the modeling efforts currently chasing the mechanistic representation of \asyn{} aggregation {\it in vivo}. To this end, both groups of models would benefit from calibration and validation against measurements of aggregate levels in living systems.
However, the availability of this quantitative information is limited by technical difficulties in monitoring the time evolution of aggregates and addressing the heterogeneous and elusive nature of oligomers \cite{Sinnige2022}.
As a result, chemical kinetic models have been trained mainly on {\it in vitro} data.
The main limitation in using this type of data relates to the unclear relevance of the results on a disease time scale \cite{Bakshi2019}, moving from the few hours or days of {\it in vitro} aggregation to the decades of PD progression. Still, studying aggregation mechanisms {\it in vitro} can provide insights into the {\it in vivo} settings. According to a recent perspective paper \cite{Sinnige2022}, the late onset of amyloid formation in sporadic PD may emerge from the intrinsic properties of aggregation kinetics, specifically from a combination of dominant secondary processes and prion-like spreading. Age-related defects in protein degradation may also work in concert with secondary events in the aggregation process to drive amyloid formation. Qualitative similarities between {\it in vitro} and {\it in vivo} aggregation mechanisms support this scenario. Indeed, biophysical features of amyloid formation detected in test-tube experiments seem to be preserved in living systems \cite{Sinnige2022, SinnigeDobson2020}. On the other hand, a quantitative analysis of the aggregation kinetics {\it in vivo} is lacking. Current studies are focusing on the applicability of {\it in vitro} molecular-level mechanistic studies to living systems and the related modeling limitations and opportunities, as extensively discussed in Meisl {\it et al.} \cite{Meisl2022}.
Molecular models combining the kinetic of \asyn{} aggregation and multiple biological pathways underlying PD neurodegeneration can assist {\it in vitro} to {\it in vivo} translation \cite{Cassotta2022}.

{\it Modulating the protein aggregation pathway}.
PD research efforts seek to provide a solid foundation for therapeutic strategy design in PD.
Not only can mechanistic models of \asyn{} homeostasis identify candidate target mechanisms, but they can also support the development of various anti-aggregation treatments currently under study. In this context, chemical kinetic models provide a theoretical framework for analyzing the inhibitory effect of multiple compounds on the aggregation process. Employing a nucleation-polymerization model that accounts for interactions with a new general component, Michaels {\it et al.} \cite{Michaels2020} identify different regimes of inhibition related to specific combinations of kinetic rate constants and binding affinities, thus optimizing the inhibitor efficacy of the compound. Furthermore, a follow-up study simulates the impact of therapeutic interventions on key oligomeric features by tuning specific rate constants \cite{Michaels2022}. This analysis suggests that inhibiting the primary nucleation event that dominates lipid-induced \asyn{} aggregation strongly reduces the peak and overall concentration of the oligomeric population. As shown in \cite{Perni2018}, this mechanism of action characterizes the natural aminosterol compound known as trodusquemine, which interferes with \asyn{} aggregation by blocking its interaction with lipid membranes and thus inhibit lipid-induced primary nucleation and fibril-mediated secondary nucleation. Also, it can suppress aggregation-related cytotoxicity by hampering oligomer interaction with cell membranes \cite{Limbocker2020}.
In addition to small molecules, potential therapeutic strategies targeting different \asyn{} species and reactions include peptides and peptidomimetics, antibodies, and molecular chaperones \cite{Fields2019}. 
In this respect, chemical kinetic studies have recently investigated the mechanism of action of specific molecular chaperones \cite{Arosio2016, Cohen2015} and antibodies \cite{Linse2020} on amyloid beta aggregation in Alzheimer's disease ({\it e.g.}, Brichos and aducanumab, respectively); the results were in good agreement with clinical outcomes. These studies should be extended to Parkinson’s disease, {\it e.g.}, supporting clinical trials that involve anti-\asyn{} antibodies such as prasinezumab and cinpanemab. 

{\it Targeting protein degradation pathways}.
Inhibiting specific molecular steps of the aggregation pathway is not the only possible intervention targeting impaired \asyn{} homeostasis. Genetic-based strategies reducing protein synthesis have been proposed to avoid aggregation in the first place, such as employing antisense oligonucleotides, RNA interference, and $\beta2$-adrenoreceptor agonists \cite{Fields2019}. 
However, side effects may arise from protein downregulation due to \asyn{} physiological role in presynaptic vesicle trafficking. Mechanistic models should be employed to determine a critical range for \asyn{} expression levels, facilitating the transition from preclinical to clinical trials. However, currently available models do not include a genetic level representation; thus, transcription and translation inhibitions can be simulated only by varying the rate constant associated with \asyn{} production. 
Another strategy consists in modulating UPP and ALP degradation mechanisms to hamper \asyn{} aggregation and PD neurodegeneration.
Autophagy can be enhanced selectively or as a whole. An overall ALP enhancement driven by chemical compounds such as rapamycin and trehalose has limited applicability in PD since it may interfere with essential physiological processes and the homeostasis of proteins other than \asyn{} \cite{Scrivo2018}. The same limitations hold for therapeutic strategies focusing on overall UPP degradation \cite{Dantuma2014}. In contrast, selective targeting of specific molecular steps in proteasomal degradation, CMA, and macroautophagy may be a promising strategy \cite{Scrivo2018}. UPP-related candidate treatments include stimulating ubiquitination and regulating the activity of chaperone systems. In addition, selective ALP activators may act by modulating LAMP2A assembly and translocation, enhancing autolysosome fusion, and targeting lysosomal function. Ideally, mechanistic models can support the development of current and future interventions of this type. However, further research is needed to reach this goal.
Most degradation models analyzed here lack a granular representation of autophagic pathways that would allow the simulation of these selective strategies. 
Still, they can provide insightful information on the role of UPP and ALP degradation pathways in neurodegeneration as the first step toward {\it in silico}-aided therapeutic strategies against PD. 

{\it Toward multiscale QSP models for PD}.
The points discussed above highlight the need for multiscale QSP models able to link molecular mechanisms to cellular and tissue-level events and, finally, clinical manifestations. These models would leverage the large amount of available data associated with neuroimaging techniques and disease-rating scales for motor symptoms (MDS-UPDRS) and cognitive impairment (MoCa and MMSE) \cite{Poewe2017, Roberts2016}. 
This type of data has already been employed by quantitative models that address the heterogeneity of the PD clinical picture by accounting only for physical symptoms \cite{Sarbaz2016}. In addition to quantitative tools supporting PD diagnosis, a mechanistic understanding of the molecular processes underlying PD pathogenesis and how they connect to the clinical endpoints is essential to inform drug discovery and development.
The models in this review represent a starting point in this sense. 
Being confined at the molecular level, they have to settle for surrogate endpoints such as ROS, \asyn{} aggregate, and LB levels to identify disease phenotypes. 
Moving to the cellular, tissue, and organ levels through network-based, epidemiological, and compartmental PK/PD approaches
\cite{Carbonell2018, Bloomingdale2021minimal, Bloomingdale2022PBPK, Lin2022}, these models can account for \asyn{} prion-like propagation between neurons and across multiple brain areas, using potential PD biomarkers such as \asyn{} concentrations in plasma and CSF \cite{Parnetti2019} and magnetic resonance imagining measurements. 

Overall, this review makes a case for using mechanistic models of \asyn{} homeostasis as valuable tools to uncover PD pathogenesis, suggesting their integration in a QSP framework to provide blueprints for therapeutic interventions. 

\section*{Conflict of Interest Statement}
The authors declare that the research was conducted in the absence of any commercial or financial relationships that could be construed as a potential conflict of interest.

\section*{Author Contributions}
ER and FR designed, wrote, and edited the manuscript. AA provided a preliminary version of some sections. LM and ED discussed the results and edited the manuscript. 
All the authors reviewed the manuscript.

\section*{Funding}
This research received no external funding.

\section*{Acknowledgments}
The authors would like to thank Dr. Gianluca Selvaggio for his helpful suggestions, comments, and support in preparing this manuscript.

\nolinenumbers

\small
\bibliography{library}

\bibliographystyle{unsrt}

\end{document}